\documentclass{amsart}
\usepackage{amscd}

\def\Bf{\mathbb}
\def\BZ{{\Bf Z}}
\def\BC{{\Bf C}}
\def\BR{{\Bf R}}

\def\zz{{\BZ/2}}

\def\zgm{{\itshape{$\BZ$-GM}}}
\def\zzgm{{\itshape{GM}}}

\def\yc{Y(\BC)}
\def\yr{Y(\BR)}
\def\xc{X(\BC)}
\def\xr{X(\BR)}
\def\H{H}
\def\alg{\mathrm{alg}}
\def\codim{\operatorname{codim}}
\newcommand{\even}{\mathrm{even}}
\newcommand{\odd}{\mathrm{odd}}
\newcommand{\stclass}{w_1(\yr)}
\newcommand{\im}{\operatorname{Im}}
\renewcommand{\ker}{\operatorname{Ker}}
\newcommand{\Coker}{\operatorname{Coker}}
\newcommand{\Gdeg}{\operatorname{deg_G}}
\newcommand{\varbar}{-}
\newcommand{\cF}{\mathcal{F}}

\newcommand{\cCinf}{\mathcal{C}^\infty}
\newcommand{\cO}{\mathcal{O}}

\newcommand{\Hom}{\operatorname{Hom}}
\newcommand{\id}{\operatorname{id}}

\newcommand{\chow}{{CH}}
\newcommand{\fcolon}{\colon\:} % the command giving the colon in f : X -> Y
\newcommand{\mapping}[5]{\begin{eqnarray*}
{#1} \fcolon {#2} & \to & {#3} \\
{#4} & \mapsto & {#5}
\end{eqnarray*}}
\newcommand{\anmapping}[4]{\begin{eqnarray*} 
{#1} & \to & {#2} \\
{#3} & \mapsto & {#4}
\end{eqnarray*}}
\def\pairing#1:#2#3#4->#5,#6#7#8|->#9{\begin{array}{rclcc}
#1 \fcolon #2 & #3 & #4 & \to & #5 \\
#6 & #7 & #8 & \mapsto & #9
\end{array}}
\def\anpairing#1#2#3->#4,#5#6#7|->#8{\begin{array}{rclcc}
#1 & #2 & #3 & \to & #4 \\
#5 & #6 & #7 & \mapsto & #8
\end{array}}

\newcommand{\iso}{\simeq}
\newcommand{\isoto}{\overset{\sim}{\to}}
\newcommand{\inject}{\hookrightarrow}

\newcommand{\labelto}[1]{\xrightarrow{\makebox[1em]{\scriptsize ${#1}$}}}
\newcommand{\Tor}{\operatorname{Tor}}
\newcommand{\Br}{\operatorname{Br}}
\newcommand{\Pic}{\operatorname{Pic}}
\newcommand{\etale}{{\text{\textrm{\'et}}}}
\newcommand{\BG}{{\Bf{G}}}
\newcommand{\BZn}{{\BZ/n}}
\newcommand{\MU}{{\boldsymbol{\mu}}}
\newcommand{\point}{\mathrm{pt}}
\DeclareMathOperator*{\tensor}{\otimes}

\newtheorem{theo}{Theorem}[section]

\newtheorem{prop}[theo]{Proposition}
\newtheorem{lem}[theo]{Lemma}
\newtheorem{cor}[theo]{Corollary}

\theoremstyle{remark}

\newtheorem{remarks}[theo]{Remarks}

\newenvironment{acknowledgements}%
{\subsubsection*{Acknowledgements}}{\par}

\begin{document}

\title[Real Enriques surfaces]{Algebraic cycles and topology \\
of real Enriques surfaces}

\author{Fr\'ed\'eric Mangolte\and Joost van Hamel}



\address{Fr\'ed\'eric Mangolte, D\'epartement de Math\'ematiques,
Universit\'e Montpellier II, 34095 Montpellier Cedex 5, France, Tel:
(33) 67 14 35 05, Fax: (33) 67 14 35 58}
\email{mangolte@math.univ-montp2.fr}

\address{Joost van Hamel, Faculteit der Wiskunde en Informatica, 
Vrije Universiteit, 
De Boelelaan 1081a, 1081 HV Amsterdam, The Netherlands,
Tel: (31) 20 444 76 94, Fax: (31) 20 444 76 53} 
\email{jvh@cs.vu.nl}

\thanks{The second author was supported in part by 
EC grant~CHRX-CT94-O506.}
 
\begin{abstract}
For a real Enriques surface $Y$ we prove that 
every homology class in $H_1(\yr, \zz)$ 
can be represented by a real algebraic curve 
if and only if all connected components of 
$\yr$ are orientable.
Furthermore, we give a characterization of 
real Enriques surfaces which are
Galois-Maximal and/or $\BZ$-Galois-Maximal and 
we determine the Brauer group of
any real Enriques surface.

\end{abstract}
\maketitle

%%%%%%%%%%%%%%%%%%%%%%%%%%%%%%%%%%%%%%%%%%%%%%%%%%%%%%%%%%%%%%%%%%%%%%
\section{Introduction}\label{sec intro}
%%%%%%%%%%%%%%%%%%%%%%%%%%%%%%%%%%%%%%%%%%%%%%%%%%%%%%%%%%%%%%%%%%%%%%

Let $Y$ be a complex algebraic surface.
Let us denote by $\yc$ the set of closed points of $Y$ endowed with
the Euclidean topology and let
$\H_2^{\alg}(\yc,\BZ)$ be the subgroup of the homology group $\H_2(\yc,\BZ)$
generated by the fundamental classes of algebraic curves on $Y$.
If $Y$ is an Enriques surface, we have
$$
\H_2^{\alg}(\yc,\BZ)= \H_2(\yc,\BZ).
$$
One of the goals of the present paper is to prove a similar property 
for real Enriques surfaces with orientable real part. 
See Theorem~\ref{theo main} below.

By an \emph{algebraic variety $Y$ over $\BR$} we mean a geometrically
integral scheme of finite type over the real numbers.  The Galois
group $G = \{1, \sigma\}$ of $\BC / \BR$ acts on $\yc$, the set of
complex points of $Y$, via an antiholomorphic involution, and the real
part $\yr$ is precisely the set of fixed points under this action.  An
algebraic variety $Y$ over $\BR$ will be called a \emph{real Enriques
surface}, a real K3-surface, etc., if the complexification $Y_\BC = Y
\otimes \BC$ is a complex Enriques surface, resp.\ a complex
K3-surface, etc.
Consider the following two classification problems:

-- classification of topological types of algebraic varieties $Y$ over
$\BR$ (the manifolds $\yc$ up to equivariant diffeomorphism),

-- classification of  topological types of the real parts $\yr$.  

For real Enriques surfaces the two classifications have been
investigated recently by Nikulin in \cite{Ni}. The topological
classification of the real parts was completed by Degtyarev and
Kharlamov who give in \cite{DeKha1} a description of all 87
topological types. Let us mention here that the real part of a real
Enriques surface $Y$ need not be connected and that a connected
component $V$ of $\yr$ is either a nonorientable surface of genus
$\leq 11$ or it is homeomorphic to a sphere or to a torus.

The problem of classifying $\yc$ up to 
equivariant diffeomorphism still lacks a satisfactory solution.
In the attempts to solve this problem, 
equivariant (co)homology
plays an important role (see \cite{Ni},  \cite{N-S}, \cite{DeKha2}).
It establishes for any algebraic variety $Y$ over $\BR$ a link between
the action of $G$ on the (co)homology of $\yc$ and the topology of $\yr$.
For example, the famous inequalities
\begin{align}
\label{eq GM1}
\dim \H_*(\yr,\zz) & \leq 
\sum_{r=0}^{2n} \dim \H^1(G, \H_r(\yc,\zz))
\\
\label{eq GM2}
\dim\H_\even(\yr,\zz) & \leq 
\sum_{r=0}^{2n} \dim \H^2(G, \H_r(\yc,\BZ)) \\  
\label{eq GM3}
\dim\H_\odd(\yr,\zz) & \leq 
\sum_{r=0}^{2n} \dim \H^1(G, \H_r(\yc,\BZ)) 
\end{align}
(cf.\ \cite{Kr1} or \cite{Si}) can be proven using equivariant homology.

We will say that $Y$ is
\emph{Galois-Maximal} or 
a \emph{\zzgm-variety} if
the first inequality turns into equality, 
and $Y$ 
will be called \emph{$\BZ$-Galois-Maximal}, or a \emph{\zgm-variety}  if
inequalities \eqref{eq GM2} and \eqref{eq GM3} are equalities.
When the
homology of $\yc$ is torsion free, 
the two notions coincide (see \cite[Prop.~3.6]{Kr1}).

A nonsingular projective 
surface $Y$ over $\BR$ with $\yr \neq \emptyset$ is both
\zzgm\ and \zgm\ if it is simply connected (see \cite[\S 5.3]{Kr1}).
If $H_1(\yc, \BZ) \neq 0$, as in the case of an Enriques surface,
the situation can be much more complicated. 
The necessary and sufficient conditions for a real Enriques surface $Y$ 
to be a \zzgm-variety were found in
\cite{DeKha2};
in the present paper we will give
necessary and sufficient conditions for
$Y$ to be \zgm. 
See Theorem \ref{theo galmax}. 

As far as we know, this is the first paper on real Enriques surfaces 
in which equivariant (co)homology with integral coefficients
is studied instead of coefficients in $\zz$.
We expect that the extra information that can be obtained
this way (compare for example equations \eqref{eq GM1}--\eqref{eq GM3} )
will be useful in the topological classification of
real Enriques surfaces.

In Section~\ref{sec brauer} we
demonstrate the usefulness of integral coefficients by computing
the Brauer group $\Br(Y)$ of any
real Enriques surface $Y$.
This completes the partial results on the $2$-torsion of $\Br(Y)$ obtained
in \cite{N-S} and \cite{N1}.
See Theorem~\ref{theo brauer}.

\subsection{Main results}

Let $Y$ be an algebraic variety over $\BR$.  Denote by
$\H^{\alg}_{n}(\yr,\zz)$ the subgroup of the homology group
$\H_{n}(\yr,\zz)$ generated by the fundamental classes of
$n$-dimensional Zariski-closed subsets of $\yr$, see \cite{BoHa} or \cite{BCR}.
We will say that these classes can be \emph{represented by algebraic
cycles.}  The problem of determining these groups 
is still open for most
classes of surfaces.

For a real rational surface $X$ we always have
$\H_2^{\alg}(\xc,\BZ)=\H_2(\xc,\BZ)$ and 
$\H_1^{\alg}(\xr,\zz)= \H_1(\xr,\zz)$, see
\cite{Si}. For real K3-surfaces, the situation is not so rigid. In
most connected components of the moduli space of real K3-surfaces the
points corresponding to a surface $X$ with 
$\dim \H_1^\alg(\xr, \zz) \geq k$ 
form a countable union of real analytic subspaces of
codimension $k$ for any $k \leq \dim \H_1(X_0(\BR), \zz)$, where $X_0$
is any K3-surface corresponding to a point from that component. In
some components this is only true for $k < \dim \H_1(X_0(\BR), \zz)$;
these components do not contain any point corresponding to a surface
$X$ with $\H_1^\alg(\xr, \zz) = \H_1(\xr, \zz)$, see \cite{Man}. For
real Abelian surfaces the situation is similar, see \cite{Huisman:abelian}.

\begin{theo}\label{theo main}
Let $Y$ be a real Enriques surface
with $\yr\ne\emptyset$. If all connected
components of the real part
$\yr$  are orientable, then
\[\H_1^{\alg}(\yr,\zz)= \H_1(\yr,\zz).\] 
Otherwise,
$$ 
\dim \H_1^{\alg}(\yr,\zz) = \dim \H_1(\yr,\zz) - 1.
$$
\end{theo}

\noindent
See Theorem~\ref{theo c alg} for more details.

In order to state further results we should mention that the set of
connected components of the real part of a real Enriques surface $Y$
has a natural decomposition into two parts $\yr=Y_1\bigsqcup
Y_2$. Following \cite{DeKha1} we will refer to these two parts as the
two \emph{halves} of the real Enriques surface.  In \cite{N1} it is
shown that $Y$ is \zzgm\ if both halves of $\yr$ are nonempty.  It follows
from
\cite[Lem.~6.3.4]{DeKha2} that if precisely one of the halves of
$\yr$ is empty, then $Y$ is \zzgm\ if and only if $\yr$ is
nonorientable.
This result plays an important role in the proof of
many of the main results of that paper (see Section~7 of \emph{loc.\
cit.}).

In the present paper we will see
in the course of proving Theorem~\ref{theo main}
that a real Enriques surface with orientable real part
is not a \zgm-variety. In Section~\ref{sec galmax} we also tackle the
nonorientable case and combining our results 
with the results for coefficients in $\zz$
that were already known we obtain the following theorem.

\begin{theo}\label{theo galmax}
Let $Y$ be a real Enriques surface with nonempty real part.
\begin{enumerate}
\item Suppose the two halves $Y_1$ and $Y_2$ are nonempty. Then $Y$ is
\zzgm. Moreover, $Y$ is \zgm\ if and only if $\yr$ is nonorientable.

\item Suppose one of the halves $Y_1$ or $Y_2$ is empty.  Then $Y$ is
\zzgm\ if and only if $\yr$ is nonorientable.  
Moreover, $Y$ is \zgm\ if and only if $\yr$ has at least one component
of odd Euler characteristic.
\end{enumerate}
\end{theo}

\noindent
There are examples of all cases described in the
above theorem (see \cite[Fig.~1]{DeKha1}).

In Section~\ref{sec brauer} we study the Brauer group $\Br(Y)$ of
a real Enriques surface $Y$ using the fact that
$\Br(Y)$ is isomorphic to the cohomological Brauer group
$\Br'(Y) = H^2_\etale(Y, \BG_m)$,
since $Y$ is a nonsingular surface.
In \cite{N-S} Nikulin and Sujatha gave various equalities and
inequalities relating
the dimension of the 2-torsion of
$\Br(Y)$ to other topological invariants of 
a real Enriques surface $Y$. It was
shown in \cite{N1} that 
$$
\dim_\zz \Tor(2, \Br(Y)) \geq 2s -1
$$
where $s$ is the number of connected components of $\yr$, and that
equality holds if $Y$ is \zzgm.
Using the results in Section~\ref{sec galmax} on 
equivariant homology with integral coefficients 
we can compute the whole group $\Br(Y)$.

\begin{theo}\label{theo brauer}
Let $Y$ be a real Enriques surface.
Let $s$ be the number of connected components of $\yr$.
If $\yr \neq \emptyset$ is nonorientable then
$$ \Br(Y) \iso (\zz)^{2s-1}. $$
If $\yr \neq \emptyset$ is orientable then
$$ \Br(Y) \iso 
\begin{cases}
(\zz)^{2s-2} \oplus \BZ/4 & \text{if both halves are nonempty},\\ 
(\zz)^{2s}                & \text{if one half is empty.}
\end{cases}
$$
If $\yr = \emptyset$ then
$$ Br(Y) \iso \zz. $$
\end{theo}

\begin{acknowledgements}
Large parts of this paper were written during visits of the first author to 
the \emph{Vrije Universiteit, Amsterdam} and of the second author to
the \emph{Universit\'e Montpellier II}. 
We want to thank J.~Bochnak and R.~Silhol 
for the invitations, and the \emph{Thomas Stieltjes Instituut} and the
\emph{Universit\'e Montpellier II} for providing the necessary funds.
We are grateful to A.~Degtyarev and V.~Kharlamov 
for giving us preliminary versions of their papers.
\end{acknowledgements}

%%%%%%%%%%%%%%%%%%%%%%%%%%%%%%%%%%%%%%%%%%%%%%%%%%%%%%%%%%%%%%%%%%%%%%
\section{Equivariant homology and cohomology }\label{sec equiv}
%%%%%%%%%%%%%%%%%%%%%%%%%%%%%%%%%%%%%%%%%%%%%%%%%%%%%%%%%%%%%%%%%%%%%%

Since the group $G = \operatorname{Gal} (\BC/\BR)$ acts in a
natural way on the complex points of an algebraic variety $Y$ 
defined over $\BR$, 
the best homology and cohomology theories for studying the topology of
$\yc$ are
the ones that take this group action into account.
In \cite{N-S} \'etale cohomology $H^*_\etale (Y, \zz)$
is used, and in \cite{N1} the observation is made that this
is essentially
the same as equivariant cohomology $H^*(\yc; G, \zz)$.
In \cite{DeKha2} Degtyarev and Kharlamov do not use 
equivariant cohomology as such, but instead a `stabilized' form
of the
Hochschild-Serre spectral sequence
$
E_{p,q}^2(X; G, \zz) = \H^p(G, H^q(X, \zz))
$.
This construction,  due to I.~Kalinin, 
is based on the fact that if $G = \zz$ then
$H^{p+2}(G, M) = H^p(G, M)$ for any group $M$ and any $p > 0$, 
and if $M$ is a $\zz$-module then even
$H^{p+1}(G, M) = H^p(G, M)$ for any $p > 0$, 
so it is possible to squeeze the
Hochschild-Serre spectral sequence into 1, or at most 2 diagonals.
They also use the analogue of this Kalinin spectral sequence in homology.
In the present paper
we stick to the original equivariant cohomology supplemented
with a straightforward dual 
construction which we call equivariant Borel-Moore homology.

First we will recall some properties of 
equivariant cohomology for a space with
an action of $G = \zz$. Then we will give the definition
of equivariant Borel-Moore homology and list the properties that we 
are going to need.
In Section~\ref{sec fund} we give a short treatment of
the fundamental class of $G$-manifolds
and formulate Poincar\'e duality in the equivariant context.

Let $X$ be a topological space with an action of $G = \zz$.
We denote the fixed point set of $X$ by $X^G$.
In \cite{Grothendieck} the groups $\H^*(X; G, \cF )$ 
are defined for a $G$-sheaf $\cF$
on $X$, which is a sheaf with a $G$-action compatible 
with the $G$-action on $X$.
Writing $G = \{1, \sigma \}$, this just means that we are given an
isomorphism of sheaves
$\varsigma:\cF \to \sigma^* \cF$
satisfying $\sigma^*(\varsigma) \circ \varsigma = \id$.
Now define
$$
\H^p(X; G, \varbar ) = R^p \Gamma(X, \varbar)^G
$$
the $p$-th right derived functor of the $G$-invariant global sections functor.
We have natural mappings 
$$
e^p_{\cF} \fcolon \H^p(X; G, \cF) \to \H^p(X, \cF)^G
$$
which are the edge morphisms of the 
\emph{Hochschild-Serre spectral sequence}
$$
E_{p,q}^2(X; G, \cF) = \H^p(G, H^q(X, \cF))
\Rightarrow \H^{p+q}(X; G, \cF)
$$

For us, the most important  $G$-sheaves will be the
constant sheaf $\zz$ and the constant sheaves 
constructed from the
 $G$-modules $\BZ(k)$ for $k \in \BZ$. 
Here we define $\BZ(k)$, to be the group of integers, equipped 
with an action of $G$
defined by
$\sigma \cdot z = (-1)^k z$.
We will use the notation $A(k)$ to denote either $\zz$ or $\BZ(k)$,
and we will sometimes use $A$ instead of $A(k)$ if $k$ is even.

There is a cup-product
$$
\H^p(X; G, A(k)) \otimes \H^q(X; G, A(l)) \to \H^{p+q}(X; G, A(k+l))
$$
and a pull-back 
$f^*$ for any continuous equivariant mapping $f \fcolon X \to Y$,
which both have the usual properties.

If $X$ is a point, 
$\H^p(\point; G, M ) = \H^p(G, M), $
which is cohomology of the group $G$ with coefficients in $M$.
Recall that
as a graded ring, $\H^*(G,\zz)$ 
is isomorphic to the polynomial ring
$\zz[\eta]$, where $\eta$ is the nontrivial element in 
$H^1(G,\zz)$. 
By abuse of notation, we will also use the notation $\eta$ for
the nontrivial element in $H^1(G, \BZ(1)) \iso \zz$ and 
$\eta^2$ for the
nontrivial element in $\H^2(G, \BZ) \iso \zz$. 
This notation
is justified by the fact that
$\eta \in H^1(G, \BZ(1))$ maps to $\eta \in H^1(G, \zz)$ 
under the reduction modulo 2 mapping and
$\eta^2 \in H^2(G, \BZ)$ maps to $\eta^2 \in H^2(G, \zz)$.

The constant mapping $X \to \point$ induces a mapping
$H^*(G, \zz) \to H^*(X; G, \zz)$ and 
we have a natural injection
$\H^p(X^G, \zz) \inject \H^p(X^G; G, \zz)$, so cup-product gives us
for any $G$-space $X$
a mapping
%\begin{equation}
  %\label{eq cdecomp}
$$ 
{\H^*(X^G,\zz)}{\otimes}{\H^*(G, \zz)}\to{\H^*(X^G;G,\zz)}
$$
%\end{equation}
which is well-known to be an isomorphism.
Taking the inverse of this isomorphism and
sending $\eta$ to the unit element in $H^*(X^G, \zz)$
we obtain a surjective homomorphism of rings
$\H^*(X^G;G,\zz) \to H^*(X^G, \zz)$
and we define for $A = \BZ$ or $\zz$ and any $k \in \BZ$
 the homomorphism of rings
$$ \beta \fcolon H^*(X; G, A(k)) \to H^*(X^G, \zz)$$
to be the composite mapping
$$ H^*(X; G, A(k)) \labelto{i^*} 
H^*(X^G; G, A(k)) \labelto{\bmod 2} H^*(X^G; G, \zz) \to H^*(X^G, \zz), $$
where $i^*$ is induced by the inclusion $i \fcolon X^G \inject X$.
Note that $\beta$ coincides with the mapping
$\beta'$ in \cite{Kr3}.
It is clear from the definition that
%\begin{equation}
  %\label{f-and-beta}
$$
  \beta(f^* \omega) = f^* \beta(\omega).
$$%\end{equation}
We use the notation 
$$
\beta^{n,p}\fcolon \H^n(X; G, A(k)) \to \H^p(X^G; \zz)
$$
for the mapping induced by $\beta$.

In Section~\ref{sec galmax}, we will need one technical lemma
which can easily be proven using the
Hochschild-Serre spectral sequence.

\begin{lem}
  \label{lem e2-not-surj}
  Let $X$ be a $G$-space with $X^G \ne \emptyset$. Then if
  $e^2_{A(k)}$ is not surjective on $\H^2(X,A(k))^G$, there is a class 
  $\omega \in \H^1(X; G, A(k-1))$ such
  that $e^1_{A(k-1)}(\omega) \ne 0$, but $\beta(\omega) = 0$. 
\end{lem}

The homology theory we are going to use is
the natural dual to equivariant cohomology.
For an extensive treatment of its properties, see \cite{JvH}.
Here we will give a short account
without proofs.

In the rest of this section we assume $X$ to be a 
locally compact space of finite cohomological dimension
with an action of $G = \zz$, and $A(k)$ will be as above.
We define the \emph{equivariant Borel-Moore homology
of $X$ with coefficients
in $A(k)$}  by
$$
H_p(X; G, A(k)) = R^{-p} \Hom_G ( R \Gamma_c (X, \BZ), A(k)) 
\text{ for $p \in \BZ$} 
$$
where $\Hom_G$ stands for homomorphisms in the category of $G$-modules
and $\Gamma_c$ stands for global sections with compact support;
this is the natural equivariant generalization of 
the usual Borel-Moore homology
in the context of sheaf theory (see, for example, \cite[Ch.IX]{Iversen}).

If $X$ is homeomorphic to an $n$-dimensional locally finite
simplicial complex
with a (simplicial) action of $G$,
the we can determine $H_p(X; G, A(k))$ from a double complex
analogous to the double complex~(1-12) in {N1}, which is used for
the calculation of equivariant cohomology.
Consider
the oriented chain complex 
$ \cCinf_n \to \cCinf_{n-1} \to \dots \to \cCinf_0 $
with closed supports 
(i.e., the elements of $\cCinf_p$ are $p$-chains that can 
be infinite).
The chain complex with coefficients in $A(k)$ is
defined by
$$ \cCinf_p(A(k)) = \cCinf_p \tensor A(k), $$
and we give it the diagonal $G$-action.
Then $H_p(X; G, A(k))$ is naturally isomorphic to to the $(-p)$th
homology group of the total complex associated to the double complex
%\begin{equation}
$$\begin{CD}
\dots  & & \dots & & \dots \\
@AAA @AAA @AAA\\
\cCinf_{n-1}(A(k)) @>{1-\sigma}>> 
\cCinf_{n-1}(A(k)) @>{1+\sigma}>> 
\cCinf_{n-1}(A(k)) @>{1-\sigma}>> 
\cdots \\
@AAA @AAA @AAA\\
\cCinf_{n}(A(k)) @>{1-\sigma}>> 
\cCinf_{n}(A(k)) @>{1+\sigma}>> 
\cCinf_{n}(A(k)) @>{1-\sigma}>> 
\cdots \\
\end{CD}$$
%\end{equation}
where the lower left hand corner has bidegree $(-n,0)$.
Note that by construction $H_p(\point; G, A(k)) = H^{-p}(G, A(k))$, 
so Poincar\'e duality holds trivially when $X$ is a point
(and the proof of Poincar\'e duality in higher dimensions, as stated
in Proposition~\ref{prop poincare}, is no more difficult than in the
nonequivariant case).
In particular,
$H_p(X; G, A(k))$ need not be zero for $p < 0$.

The groups $H_p(X; G, A(k))$ are covariantly functorial in $X$ with respect to
equivariant proper mappings and
the homomorphisms $\BZ(k) \to \zz$ induce homomorphisms
$H_p(X; G, \BZ(k)) \to H_p(X; G, \zz)$ that fit into a long exact sequence
\begin{multline}
\label{les coeff}
 \cdots \labelto{} H_p(X; G, \BZ(k)) \labelto{\times 2}
H_p(X; G, \BZ(k)) \labelto{} 
\\ \labelto{}
H_p(X; G, \zz) \labelto{}
H_{p-1}(X; G, \BZ(k)) \labelto{} \cdots
\end{multline}

As in the case of cohomology, 
there are natural homomorphisms
$$
e^{A(k)}_p \fcolon \H_p(X; G, A(k)) \to \H_p(X, A(k))^G 
$$
which are the edge morphisms
of a Hochschild-Serre spectral sequence
\begin{equation*}
  \label{spectral}
E^2_{p,q}(X; G, A(k)) = H^{-p}(G, H_q (X, A(k))) 
\Rightarrow \H_{p+q}(X; G, A(k)). 
\end{equation*}
If no confusion is likely, we use $e$ instead of $e^{A(k)}_p$;
otherwise 
we will often write
$e^{+}_p = e^{\BZ(2k)}_p$,
$e^{-}_p = e^{\BZ(2k+1)}_p$, and $e_p = e^{\zz}_p$, and we 
have similar conventions
for the edge morphisms $e^p_{A(k)}$ in cohomology.

There is a cap-product between homology and cohomology
\begin{equation*}
% \label{capproduct}
  \anpairing{\H_p(X; G, A(k))}{\otimes}{\H^q(X; G, A(l))}->%
                                     {\H_{p-q}(X; G, A(k-l))},%
   {\gamma}{\otimes}{\omega}|->{\gamma \cap \omega},
\end{equation*}
and of course we have

\begin{align}
\label{cap-and-cup}
  \gamma \cap (\omega \cup \omega') &= (\gamma \cap \omega) \cap \omega',\\
\label{cap-and-edge}
  e(\gamma \cap \omega) &= 
                    e(\gamma) \cap e(\omega),\\
%  e^{A(k-l)}_{p-q}(\gamma \cap \omega) &= 
%                    e^{A(k)}_p(\gamma) \cap e_{A(l)}^q(\omega),\\
\intertext{and for any proper equivariant mapping $f\fcolon X \to Y$}
  \label{cap-and-f}
  (f_* \gamma) \cap \omega &= f_* (\gamma \cap f^* \omega).
\end{align}

Recall that $\eta$ is the nontrivial element in $H^1(G, A(1))$.
Cap-product with $\eta$ considered as an element of $H^1(X; G, A(1))$
defines a map 
\mapping{ s^{A(k)}_p}{H_p(X; G, A(k))}{H_{p-1}(X; G, A(k+1))}%
{\gamma}{\gamma \cap \eta}

It can be shown, that the $e^{A(k)}_p$ and $s^{A(k)}_p$
fit into a long exact sequence
\begin{multline}\label{se edge}
\cdots \xrightarrow{s^{A(k-1)}_{p+1}} 
\H_p(X;G,A(k)) \xrightarrow{e^{A(k)}_p}
\H_p(X,A) \to 
\\ \to \H_p(X ;G,A(k-1))\xrightarrow{s^{A(k-1)}_{p}}
\H_{p-1}(X ;G, A(k)) \to \cdots 
\end{multline}
For $s^{A(k)}_p$ we adopt the same notational conventions as for
$e_p^{A(k)}$.

The natural mapping 
$
\H_p(X^G, A) \to \H_p(X^G; G, A)
$
and the cap-product give us a homomorphism
\begin{equation*}
  \H_*(X^G,\zz){\otimes}{\H^*(G, \zz)}\to{\H_*(X^G;G,\zz)},
\end{equation*}
which is an isomorphism.
Taking the inverse of this isomorphism and sending the nontrivial element 
$\eta \in \H^1(G,\zz)$ to the unit element in $\H^*(X^G, \zz)$
we obtain a surjective homomorphism
$$\H_*(X^G;G,\zz) \to \H_*(X^G,\zz).$$
Furthermore, the 
mapping $i_* \fcolon \H_n(X^G; G, \zz) \to \H_n(X; G, \zz)$
induced by the inclusion $i \fcolon X^G \to X$ 
is an isomorphism for any $n < 0$, so we can define
a homomorphism
$$
\rho \fcolon \H_*(X;G,A(k)) \to \H_*(X^G, \zz) 
$$
by taking the composite mapping
\begin{multline*}
  \H_*(X;G,A(k)) \labelto{\bmod 2}
  \H_*(X; G, \zz) \labelto{\cap \eta^{N}} 
  \H_{< 0} (X;G,\zz)
  \labelto{(i_*)^{-1}} 
  \\ \to 
  H_*(X^G; G, \zz) \to H_*(X^G, \zz),
\end{multline*}
where $N$ is any integer greater than the (cohomological) dimension of $X$.
We use the notation $\rho_n$ for the restriction of $\rho$ to
$\H_n(X;G,A(k))$, we write
$\rho_{n,p}$ for the composition of $\rho_n$ with
the projection $\H_*(X^G, \zz) \to \H_p(X^G, \zz)$,
and similar definitions hold for $\rho_{n, \even}$ and $\rho_{n, \odd}$.

It is clear from the above that
\begin{equation}
  \label{s-and-rho}
  \rho \circ s = \rho,
\end{equation}
and that
the mapping
$$ 
\rho_n \fcolon \H_{n}(X;G,\zz) \to \H_*(X^G,\zz)
$$
induced by $\rho$ is surjective if $n < 0$.
Note that, together with the Hochschild-Serre spectral sequence 
$E^r_{p,q}(X; G, \zz)$, this proves
equation (\ref{eq GM1}).
Equations \eqref{eq GM2} and \eqref{eq GM3} 
can be derived 
from the Hochschild-Serre spectral sequence with coefficients in $\BZ$
and the following proposition.

\begin{prop}
\label{prop zhomdecompo}
Let $X$ be a locally compact space of finite cohomological dimension
with an action of $G=\zz$.
Then
$$ 
\rho_{n,\even} \fcolon 
\H_n(X; G, \BZ(k)) \to \H_\even(X^G, \zz)
$$
is an isomorphism if $n < 0$ and $n + k$ is even, and
$$ 
\rho_{n,\odd} \fcolon 
\H_n(X; G, \BZ(k)) \to \H_\odd(X^G, \zz)
$$
is an isomorphism if $n < 0$ and $n + k$ is odd.
\end{prop}

Observe that it is not claimed that 
$\rho_n \left( \H_n(X; G, \BZ(k))  \right) \subset  \H_*(X^G,\zz)$
is contained in $\H_\even(X^G, \zz)$ (resp.\ $\H_\odd(X^G, \zz)$).
In fact this is often not the case:
for any $\gamma \in \H_n(X; G, \BZ(k))$ there is a
$p \equiv n + k \bmod 2$ such that
\begin{equation}
\label{eq bockstein}
 \rho(\gamma) = \rho_{n,p} (\gamma) + \delta(\rho_{n,p} (\gamma))
  + \rho_{n,p - 2} (\gamma) + \delta(\rho_{n,p - 2} (\gamma)) + \dotsb,
\end{equation}
where $\delta$
is the Bockstein homomorphism
$\H_{p+1}(X^G, \zz) \to \H_{p}(X^G, \zz)$
associated to the short exact sequence 
$$0 \to \BZ/2 \to \BZ/4 \to \BZ/2 \to 0$$
(compare \cite[Th.~0.1]{Kr3}).

We will also use the symbol $\delta$ for the connecting homomorphism
$\H_{n+1}(X; G, \zz) \to \H_{n}(X; G, \BZ(k))$ 
of the long exact sequence \eqref{les coeff},
and we have
\begin{align}
  \label{rho-and-delta}
   \rho_{n, \even}(\delta(\gamma))  &= 
   \rho_{n+1, \even}(\gamma) + \delta(\rho_{n+1, \odd}(\gamma)) &&
   \text{if $n + k$ is even,}\\
   \rho_{n, \odd}(\delta(\gamma))  &=
   \rho_{n+1, \odd}(\gamma) + \delta(\rho_{n+1, \even}(\gamma)) &&
   \text{if $n + k$ is odd.}
\end{align}

It is clear from the definition and the projection formula 
\eqref{cap-and-f} that
\begin{align}
   \label{cap-and-rho}
  \rho(\gamma) \cap \beta(\omega)&= \rho(\gamma \cap \omega), \\
\intertext{and for any proper mapping $f \fcolon X \to Y$ of $G$-spaces}
  \label{f-and-rho}
  \rho(f_* \gamma)&= f_* \rho(\gamma).
\end{align}

There are canonical isomorphisms
$\H_0(\point; G, A) \iso A$
and $\H_0(\point, A) = A$,
so 
the homomorphisms induced
by the constant mapping $\varphi \fcolon X \to \point$
give us for every compact $G$-space $X$  the \emph{degree maps}
\begin{align*}
  \Gdeg \fcolon \H_0(X;G, A) &{}\to A\\
\intertext{and}
  \deg \fcolon \H_0(X, A)   &{}\to A,
\end{align*}
which satisfy the equality
\begin{equation}
  \label{deg-and-edge1}
  e \circ \Gdeg = \deg{} \circ e.
\end{equation}
Extending the degree map on $\H_0(X^G, \zz)$ by $0$ to the whole of
$H_*(X^G, \zz)$,
we have by equation~(\ref{f-and-rho})
that
\begin{equation}
  \label{deg-and-rho1}
  \Gdeg(\gamma) \equiv \deg( \rho(\gamma)) \bmod2,
\end{equation}
for any $\gamma \in H_0(X; G, A)$.

Finally, define
$$
\H_*(X^G, A)^0 = \ker \left\{\deg \fcolon \H_*(X^G, A) \to A \right\},
$$
and
$\H_\even(X^G, \zz)^0 = \H_\even(X^G, \zz) \cap \H_*(X^G, \zz)^0$.
We will record three technical lemmas for use in Section~\ref{sec galmax}.
They can be proven by a careful inspection of either
the Hochschild-Serre spectral sequence $E_{p,q}(X; G, A(k))$ or 
the long exact sequence~(\ref{se edge}) with the appropriate coefficients.

\begin{lem}
\label{lem rho-zz}
Let $X$ be a compact connected $G$-space with $X^G \neq \emptyset$.
Then
$$
\rho_2 \fcolon \H_2(X; G, \zz) \to \H_*(X^G, \zz)^0
$$
is surjective if and only if the composite mapping
$$ 
\H_1 (X; G, \zz) \xrightarrow{e_1} 
 \H_1(X, \zz) ^G  \xrightarrow{\cup \eta^2}
 \H^2(G,\H_1(X, \zz))
$$
is zero.
\end{lem}

\begin{lem}
\label{lem even-rho-z}
Let $X$ be a compact connected $G$-space.
Then
$$
\rho_{2,\even} \fcolon \H_2(X; G, \BZ) \to \H_\even(X^G, \zz)^0
$$
is surjective if and only if the composite mapping
$$ 
\H_1 (X; G, \BZ(1)) \xrightarrow{e_1^-} 
 \H_1(X, \BZ(1)) ^G  \xrightarrow{\cup \eta^2}
 \H^2(G,\H_1(X, \BZ(1))
$$
is zero.
\end{lem}   

\begin{lem}
\label{lem odd-rho}
Let $X$ be a locally compact connected $G$-space with $X^G \neq \emptyset$. 
Then the mapping
$$
\rho_{2,\odd} \fcolon \H_2(X; G, \BZ(1)) \to \H_\odd(X^G, \zz)
$$
is surjective if and only if the composite mapping
$$
\H_1 (X; G, \BZ) \xrightarrow{e_1^+} \H_1(X, \BZ) ^G \xrightarrow{\cap
\eta^2} \H^2(G,\H_1(X, \BZ))
$$ 
is zero.
\end{lem}

%%%%%%%%%%%%%%%%%%%%%%%%%%%%%%%%%%%%%%%%%%%%%%%%%%%%%%%%%%%%%%%%%%%
\section{The fundamental class of a $G$-manifold} \label{sec fund}
%%%%%%%%%%%%%%%%%%%%%%%%%%%%%%%%%%%%%%%%%%%%%%%%%%%%%%%%%%%%%%%%%%%

Let again $A$ be $\zz$ or $\BZ$.
Let $X$ be an $A$-oriented topological manifold 
of finite dimension $d$ with an action of $G=\{1,\sigma\}$.
We will define the fundamental class
of $X$ in equivariant homology with coefficients in
$A(k)$ for $k$ even or odd.

It is well-known, that $H_d(X,  A) = A$, and the $A$-orientation
determines a generator $\mu_X$ of $H_d(X,  A)$.
Observe that we do not need to require $X$ to be compact, since we use
Borel-Moore homology.
If $G$ acts via an $A$-orientation preserving involution,
then $\mu_X \in \H_d(X, A)^G$, otherwise
$\mu_X \in \H_d(X, A(1))^G$.
By the Hochschild-Serre
spectral sequence \eqref{spectral} we have for $k\in\BZ$
an isomorphism $\H_d(X; G, A(k)) \iso  \H_d(X, A(k))^G$,
given by the edge morphisms $e^{A(k)}_d$,
so we a the fundamental class 
$$
\mu_X \in \H_d(X; G, A(k))
$$
where $k$ must have the right parity.

\begin{prop}[Poincar\'e duality]
  \label{prop poincare}
  Let $X$ be a $G$-manifold with fundamental class
  $\mu_X \in \H_d(X; G, A(k))$. Then the mapping
\anmapping{\H^i(X; G, A(l))}{\H_{d-i}(X; G, A(k-l))}{\omega}{\mu_X\cap\omega}
  is an isomorphism.
\end{prop}

Assuming that the action of $G$ is 
\emph{locally smooth} (i.e., each fixed point has a neighbourhood that
is equivariantly homeomorphic to $\BR^d$ with an orthogonal $G$-action), the
fixed point set of $X^G$ is again a topological manifold,
but it need not be $A$-orientable and it need not be equi-dimensional.
However, if $V$ is a connected component of $X^G$ and $V$ has
dimension $d_0$, then it has a fundamental
class $\mu_{V} \in \H_{d_0}(V, \zz)$, and we have that
the restriction of $\rho_{d,d_0}(\mu_X) \in \H_{d_0}(X^G, \zz)$ 
to $V$ equals $\mu_V$ (see \cite{JvH}).
If $X$ is a closed sub-$G$-manifold of a $G$-manifold $Y$, then the embedding
$j \fcolon X \to Y$ 
is proper, so
it induces a mapping in equivariant homology.
We define the class in $\H_d(Y; G, A(k))$ \emph{represented by $X$}
to be $j_* \mu_X$.

Now let $X$ be an algebraic variety defined over $\BR$.
We want to define the class in $\H_{2d}(X; G, \BZ(d))$
represented by a subvariety of dimension $d$.
As in \cite{Fulton}, we will distinguish two kinds of subvarieties, 
the \emph{geometrically irreducible subvarieties}, which are varieties
over $\BR$ themselves, 
and the \emph{geometrically reducible subvarieties}, which
are irreducible over $\BR$, but which split into two components
when tensored with $\BC$. Then the complex conjugation exchanges these
two components.

Any complex algebraic
variety $V$ of dimension $d$ has a fundamental class 
$\mu_V \in \H_{2d}(V(\BC), \BZ)$, and
the complex conjugation on $\BC^d$ 
preserves orientation if $d$ is even, and reverses orientation
if $d$ is odd. This implies that if $j \fcolon Z \hookrightarrow X$ is
the inclusion of a subvariety of dimension $d$ defined over $\BR$,
then $\mu_{Z_\BC}$ is a generator of $\H_d(Z(\BC), \BZ(d))^G$
if $Z_\BC$ is irreducible, and $\H_d(Z(\BC), \BZ(d))^G$
is generated by $\mu_{Z_1} + \mu_{Z_2}$ if $Z_\BC$ is
the union of two distinct complex varieties $Z_1$ and $Z_2$ 
of dimension $d$.
Hence we define the fundamental
class $\mu_Z \in \H_{2d}(Z(\BC); G, \BZ(d))$ of $Z$ to be the inverse image
of $\mu_{Z _\BC}$  
(resp.\ of $\mu_{Z_1} + \mu_{Z_2}$) under $e^{\BZ(d)}_{2d}$.
The class $[Z] \in \H_{2d}(X(\BC); G, \BZ(d))$ represented
by $Z$ is of course defined to be $j_* \mu_Z$.
If we use the notation
$[Z(\BR)] \in \H_d(\xr, \zz)$ for the homology class
represented by $Z(\BR)$, as defined in \cite{BoHa},
then indeed
\begin{equation}
  \label{rho-and-cycle1}
  \rho_{2d,d}([Z]) = [Z(\BR)]. 
\end{equation}

If $Z, Z'$ are subvarieties of $X$ defined over $\BR$ 
which are rationally equivalent
over $\BR$ (see \cite{Fulton} for a definition), then $[Z] = [Z']$, so 
we get for every $d \leq \dim X$
a well-defined cycle map 
$$\chow_d(X) \to \H_{2d}(\xc; G, \BZ(d))$$
from the Chow group in dimension $d$ to equivariant homology.
The image will be denoted by $\H_{2d}^\alg(\xc; G, \BZ(d))$, and we see by
equation~(\ref{rho-and-cycle1}), that
\begin{equation}
  \label{rho-and-cycle2}
  \rho_{2d,d}\left(\H_{2d}^\alg(\xc; G, \BZ(d))\right) =
  \H_d^\alg(\xr, \zz). 
\end{equation}

For $X$ nonsingular projective of dimension $n$, 
this map coincides with the composition of the mapping
$$\chow_d(X) \to \H^{2(n-d)}(\xc; G, \BZ(n-d))$$ 
as defined in \cite{Kr2}
and the Poincar\'e duality isomorphism.
As a consequence we can use the following description of 
the image of the cycle map in codimension 1, where we use the notation 
$\H^{2}_\alg(\xc; G, \BZ(1))$ 
for the image of $\chow_{n-1}(X)$ in cohomology.

\begin{prop}
  \label{prop lefschetz}
  Let $X$ be a nonsingular projective algebraic variety over $\BR$.
  Let $\cO_h$ be the sheaf of germs of holomorphic functions on
  $\xc$.
  Then $\H^2_\alg(\xc; G, \BZ(1))$ is the kernel of the composite mapping
  $$
  \H^2(\xc; G, \BZ(1)) \labelto{e^2_-} \H^2(\xc, \BZ) \labelto{} 
  \H^2(\xc, \cO_h)
  $$
\end{prop}

\begin{proof}
  This follows immediately from Proposition~1.3.1 in \cite{Kr2},
  which states that
  $\H^2_\alg(\xc; G, \BZ(1))$ is
  the image of the connecting morphism
  $$
  \H^1(\xc; G, \cO^*_h) \to \H^2(\xc; G, \BZ(1))
  $$
  in the long exact sequence induced by the
  exponential sequence of $G$-sheaves
  $$ 0 \to \BZ(1) \to \cO_h \to \cO^*_h \to 0 $$
\end{proof}

%%%%%%%%%%%%%%%%%%%%%%%%%%%%%%%%%%%%%%%%%%%%%%%%%%%%%%%%%%%%%%%%%%%%%%
\section{Algebraic cycles}\label{sec enriques}
%%%%%%%%%%%%%%%%%%%%%%%%%%%%%%%%%%%%%%%%%%%%%%%%%%%%%%%%%%%%%%%%%%%%%%

The following facts about real Enriques surfaces can be found in
\cite{Ni} or \cite{DeKha1}.  Let $Y$ be a real Enriques surface.
Let $X \to Y_\BC$ be the double covering 
of $Y_\BC$ by a complex K3-surface $X$. Since $\xc$ is simply connected,
$\xc$ is the universal covering space of $\yc$ and
$\H_1(\yc, \BZ) = \zz$. 
The complex conjugation $\sigma$ on $\yc$
can be lifted
to the covering $\xc$ in two different ways. If $\yr \neq \emptyset$ this
is easy to see; if $\yr = \emptyset$ we need to use the fact that 
a smooth manifold diffeomorphic to a K3-surface does not admit a free
$\BZ/4$-action, see \cite[p.~439]{Hitchin}. Hence we can give $X$ the
structure of a variety over $\BR$ in two different ways, which we will
denote by $X_1$ and $X_2$.
The two halves $Y_1$ and $Y_2$ of $\yr$ mentioned in the introduction 
consist 
of the components covered by  $X_1(\BR)$ and
$X_2(\BR)$, respectively.
All connected components of $X_1(\BR)$ and
$X_2(\BR)$ are orientable, as is the case for the
real part of any real K3-surface. 
If a connected component of a half $Y_i$ is
orientable, then it is covered by two components of $X_i(\BR)$, which
are interchanged by the covering transformation of $X$.  A
nonorientable component of $Y_i$ is covered by just one 
component of $X_i(\BR)$; this is the orientation covering.

Since for an Enriques surface $\H^2(\yc, \cO_h) = 0$ (see \cite[V.23]{BPV}),
we see by Proposition~\ref{prop lefschetz} and Poincar\'e duality
that $\H_2^\alg(\yc; G, \BZ(1)) = \H_2(\yc; G, \BZ(1))$, so
$\H_1^\alg(\yr, \zz)$ is the image of the mapping
$$ 
\alpha_2 =\rho_{2, 1} \fcolon \H_2(\yc; G, \BZ(1)) \to \H_1(\yr, \zz).
$$
In order to determine the image of $\alpha_2$
we will define $\alpha_n$
for any $n \in \BZ$ by
$$ 
\alpha_n = \rho_{n,1} \fcolon \H_n(\yc; G, \BZ(n-1)) \to \H_1(\yr, \zz).
$$
Observe, that $\alpha_n = \alpha_{n-1} \circ s^{+/-}_n$. 

\begin{lem}\label{lem codim}
For a real Enriques surface $Y$, the codimension of $\im \alpha_2$
in $\H_1(\yr,\zz)$ does not exceed $1$.
\end{lem}

\begin{proof}
We may assume that $\yr \ne \emptyset$.
Using the fact that $\alpha_{-1}$ is an isomorphism by 
Proposition~\ref{prop zhomdecompo}, and
both $s^{-}_0$ and $s^{+}_{1}$ are surjective by the 
long exact sequence~(\ref{se edge}), 
we see that $\alpha_1$ is surjective.
Since
$\alpha_1 = \alpha_2 \circ s^-_2$,
 it suffices to remark that if the cokernel of
$s^-_2\fcolon \H_2(\yc; G, \BZ(1)) \to \H_1(\yc; G, \BZ)$ is nonzero,
it is isomorphic to $\H_1(\yc, \BZ) = \zz$.
\end{proof}

\begin{prop}\label{prop rho}
Let $Y$ be a real Enriques surface. 
 
A class $\gamma\in\H_1(\yr,\zz)$ is contained in the image of 
$\alpha_2$ if and only
if
$$
\deg(\gamma \cap
\stclass) = 0,
$$ 
where $\stclass \in \H^1(\yr, \zz)$ is the first Stiefel-Whitney class
 of $\yr$.
\end{prop}

\begin{proof}
Again we may assume that $\yr \ne \emptyset$.
Denote by $\Omega$ the subspace of $\H_1(\yr,\zz)$  whose elements
$\gamma$ verify $\deg(\gamma\cap \stclass) = 0$. 

If $\yr$ is orientable, $\stclass=0$ and
$\Omega=\H_1(\yr,\zz)$. Furthermore, we have a surjective morphism
$$
\H_1(X_1(\BR),\zz) \oplus \H_1(X_2(\BR),\zz)\to\H_1(\yr,\zz)
$$
where the $X_1$ and $X_2$ are the two real K3-surfaces covering $Y$
(see the beginning of this section). This morphism fits
in a commutative diagram
$$
\begin{CD}
\H_2(X_1(\BC); G, \BZ(1))\oplus\H_2(X_2(\BC); G, \BZ(1)) @>>>
                                                    \H_2(\yc; G,\BZ(1))\\ 
@V{\alpha_2^{X_1}\oplus\alpha_2^{X_2}}VV @VV{\alpha_2}V \\ 
\H_1(X_1(\BR),\zz) \oplus \H_1(X_2(\BR),\zz)@>>>\H_1(\yr,\zz)
\end{CD}
$$
Here the $\alpha_n^{X_i} \fcolon \H_n(X_1(\BC); G, \BZ(n-1)) \to
\H_1(X_1(\BR),\zz)$ are defined in the
same way as $\alpha_n$.  As $\H_1(\xc,\BZ)=0$ for a real
K3-surface $X$, it follows from Lemma~\ref{lem odd-rho}, that
$\alpha_2^{X_1}$ and $\alpha_2^{X_2}$ are surjective, which implies
the surjectivity of $\alpha_2$. In other words, $\im \alpha_2=\Omega$.

Now assume that $\yr$ is nonorientable. 
Then $\stclass\ne 0$, and by nondegeneracy
of the cap-product pairing $\codim\Omega = 1$. 
First we will prove that $\im \alpha_2 \subset \Omega$.

Let $K=-cw_1(\yc)\in\H^2(\yc; G, \BZ(1))$, where $cw_1(\yc)$ is the first 
mixed characteristic class of the tangent bundle of $\yc$ as defined
in \cite[3.2]{Kr2}. Then $e(K) \in  \H^2(\yc,\BZ)$ is the
first Chern class of the canonical line bundle of $Y$, so $2 e(K)=0$
(see \cite[V.32]{BPV}).
This means that
for any $\gamma\in\H_2(\yc; G, \BZ(1))$ we have
$\Gdeg(\gamma\cap K)=\deg(e(\gamma) \cap e(K)) = 0$,
so
$\deg(\rho(\gamma) \cap \beta(K)) = 0$ by equations \eqref{deg-and-rho1}
and \eqref{cap-and-rho}.

The projection $\rho_{2,2}(\gamma)$ of
$\rho(\gamma) \in \H_*(\yr, \zz)$ to $\H_2(\yr, \zz)$ is zero 
by equation \eqref{eq bockstein}
and
the projection $\beta^{2,0}(K)$ of $\beta(K) \in \H^*(\yr, \zz)$ 
to $\H^0(\yr, \zz)$ is zero by \cite[Th.~0.1]{Kr3}.
This implies
$$
\deg(\rho(\gamma) \cap \beta(K)) = \deg (\rho_{2,1}(\gamma) \cap
\beta^{2,1}(K) ),
$$ 
but $\beta^{2,1}(K) = \stclass$ by
\cite[Th.~3.2.1]{Kr2}, and $\rho_{2,1}(\gamma)= \alpha_2(\gamma)$ by
definition, so $\deg(\alpha_2(\gamma) \cap \stclass)=0$.  
In other words, $\im \alpha_2\subset\Omega$.  
Lemma~\ref{lem codim} now gives us that $\im \alpha_2 = \Omega$.
\end{proof}

\begin{cor}
  \label{cor alpha}
  With the above notation,
  $\alpha_2$ is surjective if and only if
  $Y(\BR)$ is orientable.
\end{cor}

Theorem~\ref{theo main} in the introduction is an immediate
consequence of 
Proposition~\ref{prop rho}. We can even give an 
explicit description of $\H_1^\alg(\yr, \zz)$.

\begin{theo}\label{theo c alg}
Let $Y$ be a real Enriques
surface. 

A class $\gamma\in\H_1(\yr,\zz)$
can be represented by an algebraic cycle 
if and only if 
$\deg(\gamma \cap\stclass) = 0$.
\end{theo}

%%%%%%%%%%%%%%%%%%%%%%%%%%%%%%%%%%%%%%%%%%%%%%%%%%%%%%%%%%%%%%%%%%%%%%
\section{Galois-Maximality}\label{sec galmax}
%%%%%%%%%%%%%%%%%%%%%%%%%%%%%%%%%%%%%%%%%%%%%%%%%%%%%%%%%%%%%%%%%%%%%%

The aim of this section is to describe which Enriques surfaces
are \zgm-varieties and/or \zzgm-varieties in terms of the 
orientability of the real part and the distribution of 
the components over the halves. See 
the introduction
for the definition of Galois-Maximality and
Section \ref{sec enriques} for the
definition of 'halves'.

The proof of Theorem \ref{theo galmax} will consist of a collection of
technical results and explicit constructions
of equivariant homology classes.
For completeness we also prove the 
parts of Theorem \ref{theo galmax} concerning coefficients in $\zz$, 
although these results are not new (see the introduction).

\begin{lem}\label{caract gm}
Let $Y$ be an algebraic variety over $\BR$. Then  
\begin{enumerate}
\item $Y$ is \zgm\ if and only if $e^+_p$ is surjective
on $\H_p(\yc,\BZ)^G$ and $e^-_p$ is surjective onto
$\H_p(\yc,\BZ(1))^G$ for all $p$.
\item $Y$ is \zzgm\ if and only $e_p$ is surjective onto
$\H_p(\yc,\zz)^G$ for all $p$.
\end{enumerate}
\end{lem}

\begin{proof}
  This follows from the fact that $Y$ is \zzgm\ (resp. \zgm)
  if and only if the Hochschild-Serre spectral sequence
  $E^r_{p,q}(\yc; G, A)$ is trivial for $A = \zz$ (resp. $\BZ$),
  and this can be checked by looking at the edge morphisms,
    since we have for every $k \geq 0$ and every $G$-module $M$ 
    natural surjections
    $H^k(G, M) \to H^{k+2}(G, M)$, and 
    $H^k(G, M) \to H^{k+1}(G, M(1))$, which are isomorphisms for
    $k > 0$.
\end{proof}

\begin{lem}\label{lem enr gm}
Let $Y$ be a real Enriques surface with $\yr\ne\emptyset$. Then
\begin{enumerate}
\item
for any $p\in\{0,2,3,4\}$, $e^{+/-}_p$  is surjective onto  
$\H_p(\yc,\BZ(k))^G$,

\item
for any $p\in\{0,3,4\}$, $e_p$ is surjective onto
$\H_p(\yc,\zz)^G$.
\end{enumerate}
\end{lem}

\begin{proof}  
  This can be seen from the Hochschild-Serre spectral sequences
  (cf. \cite[\S~5]{Kr1}).
\end{proof}

\begin{cor}
  \label{cor enr gm}
  Let $Y$ be a real Enriques surface with $\yr\ne\emptyset$. Then $Y$ is
  \zgm\ if and only if $e_1^{+/-}$ is surjective onto
  $\H_1(\yc, \BZ(k))^G$ for $k=0, 1$.
  Moreover, $Y$ is \zzgm\ if and only if $e_1$ and $e_2$
  are surjective onto $\H_1(\yc, \zz)^G$, resp.\ $\H_2(\yc, \zz)^G$.
\end{cor}

\begin{lem}\label{e2e1}
Let $Y$ be a real Enriques surface with $\yr\ne\emptyset$.
If $e_2$ is not surjective onto $\H_2(\yc,\zz)^G$, then
$e_1$ is not surjective onto $\H_1(\yc,\zz)^G$.
\end{lem}

\begin{proof}
By Poincar\'e duality
we see that
if $e_2$ is not surjective onto $\H_2(\yc,\zz)^G$, then
 $e^2$ is not surjective onto $\H^2(\yc,\zz)^G$. 
Let us assume that $e^2$ is not surjective.
Then by Lemma~\ref{lem e2-not-surj} there exists an
$\omega \in \H^1(\yc; G, A(k-1))$ such
that $e^1_{A(k-1)}(\omega) \ne 0$, but $\beta(\omega) = 0$. 

Now suppose $e_1$ is surjective onto $\H_1(\yc,\zz)^G$, then there exists
a $\gamma\in \H_1(\yc;G,\zz)$ such that
$$
\deg(e_1(\gamma)\cap e^1(\omega))\ne 0.
$$
This means that $\Gdeg(\gamma \cap \omega) \ne 0$, but this contradicts
$$
\Gdeg(\gamma \cap \omega) = \deg(\rho(\gamma) \cap \beta(\omega))
= \deg(\rho(\gamma) \cap 0) = 0.
$$
Hence $e_1$ is not surjective.
\end{proof}

\begin{prop}\label{prop zgmzzgm}
Let $Y$ be a real Enriques surface with $\yr\ne\emptyset$.
Then
\begin{enumerate}
\item $Y$ is \zgm\ if and only if $e^+_1$ and $e^-_1$ are nonzero.
\item $Y$ is \zzgm\ if and only if $e_1$ is nonzero.
\item If $e_1$ is zero then $e^+_1$ and $e^-_1$ are zero.
In particular, if $Y$ is \zgm, then $Y$ is also \zzgm.
\end{enumerate}
\end{prop}

\begin{proof}
If $Y$ is an Enriques surface, 
$$
\H_1(\yc,\BZ)=\H_1(\yc,\zz)=\zz,
$$
so $e^{+/-}_1$ and $e_1$ are surjective if and only if they are
nonzero. By Lemma~\ref{e2e1},  
$e_2$ is surjective if $e_1 \ne 0$, so 
we obtain the
first two assertions from Corollary~\ref{cor enr gm}. 
The last assertion follows from the commutative diagram
$$\begin{CD}
\H_1(\yc; G, \BZ(k)) @>{e^{+/-}_1}>> \H_1(\yc, \BZ(k))\\
@VVV                                    @VVV          \\
\H_1(\yc; G, \zz)    @>{e_1}>>      \H_1(\yc, \zz)
\end{CD}$$

\end{proof}

\begin{lem}\label{lem e1+}
Let $Y$ be a real Enriques surface with $\yr\ne\emptyset$. 
Then $e_1^+=0$ if and only if $\yr$ is orientable.
\end{lem}
\begin{proof}
We know from Corollary~\ref{cor alpha}, that
$\alpha_2$ is surjective if and only if
$\yr$ is orientable.
Since $\H_1(\yc,\BZ)= \zz$, the mapping
$H_1(\yc, \BZ) ^G \xrightarrow{\cup \eta^2} 
\H^2(G,\H_1(\yc, \BZ))$ is an isomorphism,
so Lemma~\ref{lem odd-rho} gives us that
$\alpha_2$ is surjective if and only if $e_1^+=0$.
\end{proof}

\begin{lem}\label{lem e1-}
If the two halves $Y_1$ and $Y_2$ of a real Enriques surface
$Y$ are nonempty,
then $e^-_1\ne 0$.
\end{lem}
\begin{proof}
Let $X$ be the K3-covering of $Y_\BC$, 
let $\tau$ be the deck transformation of this covering
and denote by
$\sigma_1$ and $\sigma_2$ the two different
involutions of $\xc$ lifting the involution
$\sigma$ of $\yc$. Let $X_i(\BR)$ be the set of fixed points under
$\sigma_i$ and let $p_i$ be a point in $X_i(\BR)$ for $i = 1,\;2$. 

Let $l$ be an arc in $\xc$ connecting $p_1$ and $p_2$ without containing
any other point of $X_1(\BR)$ or $X_2(\BR)$.
Then the union $L$ of the four arcs 
$l$, $\sigma_1(l)$, $\sigma_2(l)$, 
$\tau(l)$ is homeomorphic to a circle, and we have that $\tau(L)=L$.
This implies that the image $\lambda$
of $L$ in $\yc$
is again homeomorphic to a circle; we choose an orientation on $\lambda$.

Now $G$ acts on $\lambda$ via an orientation reversing involution,
so $\lambda$ represents a class $[\lambda]$ in $\H_1(\yc; G, \BZ(1))$.
Since $\xc \to \yc$ is the universal covering,
and the inverse image of $\lambda$ is precisely $L$, hence homeomorphic
to a circle, the class of $\lambda$ is nonzero in $\H_1(\yc,\BZ)$,
so $e^-_1([\lambda]) \ne 0$.
\end{proof}

\begin{lem}\label{lem one or}
If exactly one of the halves $Y_1$, $Y_2$ of a real Enriques surface
$Y$ is empty, then $e_1 = 0$ if and only if $\yr$ is orientable.
\end{lem}

\begin{proof}
If $e_1=0$, we have $e_1^+=0$ by Proposition~\ref{prop zgmzzgm}
and then $\yr$ is orientable by \ref{lem e1+}.  Conversely, if $\yr$ is
orientable and $X_2(\BR) = \emptyset$, then $X_1(\BR) \to Y(\BR)$ is
the trivial double covering, so it induces a surjection
$\H_*(X_1(\BR), \zz)^0 \to \H_*(\yr, \zz)^0$, where $\H_*({-}, \zz)^0$
denotes the kernel of the degree map as defined in Section~\ref{sec
equiv}.  Since $H_1(\xc, \zz) = 0$, the mapping $\rho\fcolon
\H_2(X_1(\BC); G, \zz) \to \H_*(X_1(\BR), \zz)^0$ is surjective by
Lemma~\ref{lem rho-zz}.  Now the functoriality of $\rho$ with respect
to proper equivariant mappings (equation (\ref{f-and-rho})) implies
$$
\rho_2 \fcolon\H_2(Y(\BC); G, \zz) \to \H_*(Y(\BR), \zz)
$$ 
is surjective, and Lemma~\ref{lem rho-zz} then gives that $e_1$ is
zero.
\end{proof}

\begin{lem}\label{lem euler char}
If exactly one of the halves $Y_1$, $Y_2$ of a real Enriques surface
$Y$ is empty, then $e^-_1 \ne 0$ if and only if
$\yr$ has components of odd Euler characteristic.
\end{lem}

\begin{proof}
Assume $Y_2 = \emptyset$.
By Lemma~\ref{lem even-rho-z}, it suffices to show that
$$
\rho_{2,\even} \fcolon \H_2(\yc; G, \BZ) \to \H_\even(\yr, \zz)^0
$$
is surjective if and only if
$\yr$ has no components of odd Euler characteristic.
Although $\yr$ need not be orientable, 
we can apply the K3-covering as in the previous lemma and
prove that the image of $\rho_{2,\even}$
contains a basis for the subgroup $\H_0(\yr, \zz) \cap \H_\even(\yr, \zz)^0$,
so $\rho_{2,\even}$ is surjective if and only if
$$
\rho_{2,2} \fcolon \H_2(\yc; G, \BZ) \to  \H_2(\yr, \zz)
$$
is surjective.
We will use that $\H_2(\yr, \zz)$ is generated by 
the fundamental classes of the connected components
of $\yr$.

Pick a component $V$ of $\yr$. If $V$ is orientable, it
gives a  class in $\H_2(\yc; G, \BZ)$, which maps to the fundamental class
of $V$ in  $\H_2(\yr, \zz)$.
Now assume $V$ is nonorientable.
Let $[V]$ be the fundamental class of $V$ in $H_2(\yr, \zz)$,
let $[V]_G$ be the class represented by $V$ in $H_2(\yc; G, \zz)$, and
let $\gamma = \delta([V]_G)$ be the Bockstein image in $H_1(\yc; G, \BZ(1))$.
Then $\rho_{1,2}(\gamma) = \rho_{2,2}([V]_G) = [V]$ by
equation~(\ref{rho-and-delta}), so $[V]$ is in the image of
$H_2(\yc; G, \BZ)$ under $\rho_{2,2}$ if and only if $e^-_1(\gamma) = 0$.

From the construction of $\gamma$ we see that 
$e^-_1(\gamma) = i_* \delta([V])$, where $i \fcolon V \to \yc$ 
is the inclusion and $\delta([V]) \in H_1(V, \BZ)$ is the
Bockstein image of $[V]$. Therefore  $e^-_1(\gamma)$ can be represented
by a circle $\lambda$ embedded in $V$.
Since $\xc \to \yc$ is the universal covering,
$e^-_1(\gamma)$ is zero if and only if the inverse image $L$ of $\lambda$
in $\xc$ has two connected components.
Let $W$ be the component of $X_1(\BR)$ covering $V$. Then $W$ is
the orientation covering of $V$ and $L \subset W$.
If $V$ has odd Euler characteristic, 
then it is the connected sum of a real projective plane
and an orientable compact surface. We see by elementary geometry
that $L$ is connected.
If $V$ has even Euler characteristic, it is the connected sum of a 
Klein bottle and an orientable compact surface, 
and we see that $L$
has two connected components.
\end{proof}

\begin{proof}[Proof of Theorem \ref{theo galmax}]
By Proposition~\ref{prop zgmzzgm}, the first part of the theorem
follows from Lemma~\ref{lem e1+} and Lemma~\ref{lem e1-}, 
and the second part of the theorem follows from Lemma~\ref{lem one or}
and lemma~\ref{lem euler char}.
\end{proof}

%%%%%%%%%%%%%%%%%%%%%%%%%%%%%%%%%%%%%%%%%%%%%%%%%%%%%%%%%%%%%%%%%%%%%%
\section{The Brauer group}\label{sec brauer}
%%%%%%%%%%%%%%%%%%%%%%%%%%%%%%%%%%%%%%%%%%%%%%%%%%%%%%%%%%%%%%%%%%%%%%

Let $Y$ be a nonsingular projective algebraic variety defined
over $\BR$.
Let $$\Br'(Y) = H^2_\etale(Y, \BG_m)$$
be the cohomological Brauer group of $Y$, and let $\Tor(n, \Br'(Y))$
be the $n$-torsion of $\Br'(Y)$.
We have a canonical isomorphism 
\begin{multline}
\label{eq torbr}
\Tor(n, \Br'(Y)) \iso \\
\iso \Coker \left \{H^2_\alg(\yc; G, \BZ(1)) 
\labelto{\bmod n} H^2(\yc; G, \BZn(1)) \right \},
\end{multline}
as can be deduced from the Kummer sequence
$$ 1 \labelto{} \MU_n \labelto{} 
\BG_m \labelto{\times n} \BG_m \labelto{} 1, $$
and the well-known identifications
\begin{align*}
H^k_\etale(Y, \MU_n) & \iso H^k(\yc; G, \BZn(1)) \\
\intertext{and}
H^1(Y, \BG_m) & \iso \Pic(Y).
\end{align*}
It can be checked, that the mapping
$$\beta^{2,0} \fcolon H^2(\yc; G, \zz) \to H^0(\yr, \zz)$$
induces a well-defined homomorphism
\begin{equation}
\label{CTP}
\Tor(2, \Br'(Y)) \to H^0(\yr, \zz).
\end{equation}

If $\dim Y \leq 2$, in particular if $Y$ is
a real Enriques surface, we may identify $\Br'(Y)$
with the classical Brauer group $\Br(Y)$ (see \cite[II, Th.~2.1]{Brauer}).
Two of the main problems considered in \cite{N-S} and \cite{N1} 
are the calculation
of $\dim_\zz \Tor(2, \Br(Y))$ and the question whether the mapping
\eqref{CTP} is surjective for every real Enriques surface $Y$. 
Both problems were solved for certain classes of real Enriques surfaces.
 The second problem has been completely solved
in \cite{Kr3}, where it is shown that the mapping \eqref{CTP} is
surjective for any nonsingular projective 
surface $Y$ defined over $\BR$ (see
Remark~3.3 in \emph{loc.\ cit.}).  
The results in Section~\ref{sec galmax} will help us
to solve the first problem for every Enriques
surface $Y$ by determining the whole
group $\Br(Y)$.

\begin{lem}\label{lem torbr}
Let $Y$ be a nonsingular projective algebraic variety defined over 
$\BR$ such that
$$ H^2_\alg((\yc; G, \BZ(1)) = H^2(\yc; G, \BZ(1)). $$
Then
$$ \Tor(\Br'(Y)) \iso \Tor(H^3(\yc; G, \BZ(1))). $$
\end{lem}
\begin{proof}
By the hypothesis and the isomorphism \eqref{eq torbr} there is
for every integer $n > 0$ a short exact sequence
$$ H^2(\yc; G, \BZ(1)) \tensor \BZn \to H^2(\yc; G, \BZn(1))
\to \Tor(n, \Br'(Y)), $$
hence we deduce from the long exact sequence in equivariant cohomology
associated to the short exact sequence
$$ 0 \to \BZ(1) \labelto{\times n} \BZ(1) \to \BZn(1) \to 0 $$
that we have for every $n > 0$ a natural isomorphism
$$ \Tor(n, \Br'(Y)) \iso \Tor(n, H^3(\yc; G, \BZ(1))). $$
\end{proof}

\begin{proof}[Proof of Theorem \ref{theo brauer}]
By \cite[I.2 and II, Th.~2.1]{Brauer} we have  
$\Br(Y) = \Tor(\Br(Y)) = \Tor(\Br'(Y))$.
On the other hand,
$\Tor(H^3(\yc; G, \BZ(1))) = H^3(\yc; G, \BZ(1)) $
since $H^3(\yc, \BZ) = \zz$. Hence, by Lemma~\ref{lem torbr}
and Poincar\'e duality
$$\Br(Y) \iso H_1(\yc; G, \BZ(1)). $$
Now consider the long exact sequence~\eqref{se edge} for $A(k) = \BZ$:
$$ \dots \labelto{e_1^+} H_1(\yc, \BZ) \to
H_1(\yc; G, \BZ(1))
\labelto{s_1^-} H_0(\yc; G, \BZ) \to \dotsb $$
It follows from Proposition~\ref{prop zhomdecompo} and the long
exact sequence~\eqref{se edge} for $A(k) = \BZ(1)$ that 
$\rho \fcolon H_*(\yc; G, \BZ) \to H_*(\yr, \zz)$ induces an isomorphism
$$ \im s_1^- \isoto H_\even(\yr, \zz)^0.$$
We obtain an exact sequence
\begin{equation}\label{se br}
 \dots \labelto{e_1^+} \zz \to 
H_1(\yc ; G, \BZ(1)) \to  H_\even(\yr, \zz)^0\to 0.
\end{equation}

\noindent
If $\yr \neq \emptyset$ 
is nonorientable, then $e_1^+ \neq 0$ by Lemma~\ref{lem e1+},
so $H_1(\yc ; G, \BZ(1)) \iso (\zz)^{2s-1}$, 
which proves the first part of the theorem.

Now assume $\yr \neq \emptyset$ is orientable. 
Then $e_1^+ = 0$ by Lemma~\ref{lem e1+},
so we get from~\eqref{se br} an exact sequence
\begin{equation*}
0 \to \zz \to H_1(\yc; G, \BZ(1)) \to (\zz)^{2s -1} \to 0,
\end{equation*}
hence 
$H_1(\yc; G, \BZ(1)) \iso (\zz)^{2s}\ 
\text{or}\ (\zz)^{2s-2} \oplus (\BZ/4).$

In order to decide between these two possibilities, 
consider the following commutative diagram with exact rows.
$$\begin{CD}
H_2(\yc; G, \zz) @>{\delta^-}>> 
H_1(\yc; G, \BZ(1)) @>{\times 2}>>
H_1(\yc; G, \BZ(1)) \\
@V{e_2}VV  @V{e_1^-}VV   @V{e_1^-}VV \\
H_2(\yc, \zz) @>{\delta}>> 
H_1(\yc, \BZ) @>{\times 2}>>
H_1(\yc, \BZ) \\
@A{e_2}AA  @A{e_1^+}AA  @A{e_1^+}AA \\
H_2(\yc; G, \zz) @>{\delta^+}>> 
H_1(\yc; G, \BZ) @>{\times 2}>>
H_1(\yc; G, \BZ)
\end{CD}$$
We have that $H_1(\yc; G, \BZ(1))$ is pure $2$-torsion
if and only if $\delta^-$ is surjective. We
claim that $\delta^-$ is surjective if and only if
$e_1^- = 0$.
Together with Lemmas~\ref{lem euler char} and \ref{lem e1-} 
this  would prove the second part of the theorem.

Let us prove the claim.
Since $e_1^+ = 0$, we have $\delta \circ e_2 =0$.
If $e_1^- \neq 0$, an easy diagram chase shows that 
$\delta^-$ is not surjective.
On the other hand
the following diagram can be shown to be commutative.
\[
\setlength{\unitlength}{0.00083300in}
\begin{picture}(4700,1200)(58,-319) 
\thinlines
\put(751,464){\vector( 0,-1){525}}
\put(1576,539){\vector( 3,-1){1350}}
\put(1576,-211){\vector( 1, 0){1350}}
\put(751,614){\makebox(0,0)[b]{\smash{$H_2(\yc; G, \BZ)$}}}
\put(751,-286){\makebox(0,0)[b]{\smash{$H_2(\yc ; G, \zz)$}}}
\put(3226,-286){\makebox(0,0)[lb]{\smash{$H_1(\yc; G, Z(1))$}}}
\put(2176,-136){\makebox(0,0)[b]{\smash{$\scriptstyle \delta^-$}}}
\put(2176,464){\makebox(0,0)[lb]{\smash{$\scriptstyle s^+_2$}}}
\put(676,239){\makebox(0,0)[rb]{\smash{$\scriptstyle \bmod 2$}}}
\end{picture}
\]

In other words,  $\im s_2^+ \subset \im \delta^-$.
Now $\ker e_1^- = \im s_2^+$, so if $e_1^- = 0$, then $\delta^-$ is
surjective.

Finally, we will consider the
short exact sequence~\eqref{se br}
for the 
case $\yr = \emptyset$.
Then $G$ acts freely on $\yc$, so we have for all $k$ that
$H_k(\yc; G, \zz) = H_k(\yc/G, \zz)$.
By the remarks made in the introduction
of Section~\ref{sec enriques}, this means that
$H_1(\yc; G, \zz) = \zz \times \zz$,
and we can see from the long exact sequence~\eqref{se edge} for
$A(k) = \zz$ that $e_1 = 0$. 
This implies that
$e^+_1 = 0$ (see Proposition~\ref{prop zgmzzgm}.iii), 
hence $H_1(\yc; G, \BZ(1)) = \zz$.
\end{proof}

%%%%%%%%%%%%%%%%%%%%%%%%%%%%%%%%%%%%%%%%%%%%%%%%%%%%%%%%%%%%%%%%%%%%%%

\end{document}